\documentclass[12pt]{article}

\usepackage{amssymb, amsmath, amsfonts}
\usepackage{mathrsfs}
\usepackage{multicol}
\usepackage{graphicx}
\usepackage{amsthm}
\usepackage[english]{babel}

\makeatletter
\renewenvironment{table}[1][]
  {\def\@captype{table}}
  {}

\renewenvironment{figure}
  {\def\@captype{figure}}
  {}
\makeatother

\newtheorem{df}{Definition}
\newtheorem{Lem}{Lemma}
\newtheorem{Th}{Theorem}
\newtheorem{Cor}{Corollary}

\theoremstyle{remark}
\newtheorem{Rem}{Remark}

\newcommand{\E}{{\sf E}}
\newcommand{\p}{{\sf P}}

\newcommand{\ind}{\mathbb I}

\newcommand{\Co}{{\sf cov}}

\begin{document}
\centerline{SIMULATION  AND  ANALYTICAL  APPROACH  TO THE IDENTIFICATION}

\centerline{OF SIGNIFICANT FACTORS}
\vskip 3mm

\vskip 5mm
\noindent Alexander V. Bulinski$^{1}$ and Alexander S. Rakitko
\!\!\footnote{The work is partially supported by RFBR grant 13-01-00612.}

\noindent Faculty of Mathematics and Mechanics,
Lomonosov Moscow State University

\noindent Moscow 119991, Russia

\noindent bulinski@yandex.ru

\vskip 8mm

\noindent ABSTRACT

We develop our previous works concerning the identification of the
collection of significant factors determining some, in general,
non-binary random response variable. Such identification is
important, e.g., in biological and medical studies.
Our approach is to examine the quality of response variable
prediction by functions in (certain part of) the factors.
The prediction error estimation requires some
cross-validation procedure, certain prediction algorithm and
estimation of the penalty function. Using simulated data we
demonstrate
the efficiency
of our method.
We prove a new
central limit theorem for introduced regularized estimates
under some natural conditions for arrays of exchangeable random variables.

\vskip 10mm

\noindent 1.   INTRODUCTION

In a number of models the (random) response variable $Y$
depends on some factors $X_1,\ldots,X_n$. A nontrivial problem is to
identify the set of the most ``significant factors''.
Loosely speaking, for a given
$r<n$ one can try to find such collection $\{k_1,\ldots,k_r\}
\subset \{1,\ldots,n\}$ that $Y$ depends   ``essentially'' on
$X_{k_1},\ldots,X_{k_r}$ and the impact of other factors can be
viewed as negligible.
Note that the problem of this type is important in medical and
biological studies where $Y$ can describe the state of a patient
health. For instance, $Y=1$ or $Y=-1$ may indicate that a person is
sick or healthy, respectively. Note also that in pharmacological
studies the values $-1$ or $1$ of a response variable can describe
efficient or inefficient application of some medicine. Thus it is
clear that binary response variables play an important role in
various disciplines. At the same time it is obvious that more
detailed description of experiments can be desirable. In this regard
we refer, e.g., to Bulinski and Rakitko (2014) where non-binary response
variables were studied.

There exist various complementary approaches concerning the
prediction of response variable and selection of  significant
combinations of factors. Such analysis in medical and biological
studies is included in special research domain called the {\it
genome-wide association studies} (GWAS). The problems and progress
in this important domain are considered, e.g., in Moore et al. (2010) and
Visscher et al. (2012). Among powerful statistical tools applied in GWAS
one can indicate the principle component analysis (Lee et al. (2012)),
logistic and logic regression (Schwender and Ruczinski (2010),
Sikorska et al. (2013)), LASSO (Tibshirani and Taylor (2012)) and various
methods of statistical learning (Hastie et al. (2008)). Note also that
there are various modifications of these methods.

We are  interested in the ``dimensionality reduction'' of the whole
collection of factors and so employ the term ``MDR method''. This
term was introduced, for binary response variable, in the paper
Ritchie et al. (2001)  and goes back to the Michalski
algorithm. However, instead of considering
contiguity tables (to specify zones of low and high risk) presented
in Ritchie et al. (2001) and many subsequent works we choose another way.
Namely, to predict (in general non-binary) $Y$ we use some function
$f$ in factors $X_1,\ldots,X_n$. The quality of such $f$ is
determined by means of the error function $Err(f)$ involving a
penalty function $\psi$. This penalty function allows us to take
into account the importance of different values of $Y$. As the law
of $Y$ and $X=(X_1,\ldots,X_n)$ is unknown we cannot evaluate
$Err(f)$. Thus statistical inference is based on the estimates of
error function. Developing Bulinski et al. (2012), Bulinski (2012), Bulinski (2014)  we propose (in more general setting) statistics constructed by means of a
prediction algorithm for response variable and $K$-fold
cross-validation procedure. One of the main results of
Bulinski and Rakitko (2014) gives the criterion of strong consistency of the
mentioned error function estimates when the number of observations
tends to infinity. The strong consistency is essential because to
identify the ``significant collection'' of factors we have to
compare simultaneously a number of statistics. Moreover, we proposed
in Bulinski (2014) and Bulinski and Rakitko (2014) the {\it regularized
versions} of the employed statistics (involving the appropriate
estimates of the penalty function) to establish the central limit
theorem (CLT).

The paper is organized as follows. Section 2 contains notation and
auxiliary results.  In Section 3 we discuss the results of
simulations to identify (according to our method) the collection of
significant factors determining a binary response variable. In
Section 4 we prove the new CLT for our estimates (in general for
non-binary response $Y$) using some natural conditions concerning
the arrays of exchangeable random variables.

\vskip 10mm

\noindent 2. NOTATION AND AUXILIARY RESULTS

Further on we suppose that all random variables under consideration
are defined on a probability space $(\Omega,\mathcal{F},\p)$. Let
$Y$ take values in a finite set $\mathbb{Y}$ which we will identify
with the set $\{-m,\ldots,m\}$ where $m \in \mathbb{N}$. To comprise
binary variables we can assume that their values belong to the set
$\{-1,0,1\}$ and the value $0$ is taken with probability $0$. Let
also $X_1,\ldots,X_n$ take values in  an arbitrary finite set
$\mathbb{X}=\{0,\ldots,s\}$. Choose $f:\mathbb{X}\to \mathbb{Y}$ and
a penalty function $\psi: \mathbb{Y}\to \mathbb{R}_+$. The trivial
case $\psi \equiv 0$ is excluded.
Introduce the {\it error function}
$$
Err(f):={\sf E}|Y-f(X)|\psi(Y).
$$
It is easily seen that one can write $Err(f)$ in the following way
\begin{equation*}
Err(f)=\sum_{y,z\in \mathbb Y} |y-z|\psi(y){\sf P}(Y=y,f(X)=z)=
\sum_{z\in \mathbb{Y}} \sum_{x\in A_z} w^{\top}(x)q(z)
\end{equation*}
where $q(z)$ is the $z$-th column of  $(2m+1)\times (2m+1)$ matrix
$Q$ with entries $q_{y,z}=|y-z|$, $y,z\in\mathbb{Y}$ (the entry
$q_{-m,-m}$ is located at the left upper corner of $Q$),
$$
w(x)=(\psi(-m){\sf P}(Y=-m,X=x),\ldots,\psi(m){\sf
P}(Y=m,X=x))^{\top}
$$
and $\top$ stands for transposition. All vectors are considered as
column-vectors.
According to Bulinski and Rakitko (2014) we can rewrite $Err(f)$ as follows
\begin{equation}\label{Err}
Err(f)
=\sum_{i=0}^{2m-1}\sum_{i-m< |y| \le m}\psi(y){\sf
P}(Y=y,|f(X)-y|>i).
\end{equation}
    The law of $(X,Y)$ is unknown, therefore, for each $f: \mathbb{X}\to
\mathbb{Y}$, we can not evaluate $Err(f)$. Thus it is natural that
statistical inference concerning the quality of prediction of the
response variable $Y$ by means of $f(X)$ is based on the estimates
of $Err(f)$.

Let $\xi^1,\xi^2,\ldots$ be a sequence of independent identically
distributed (i.i.d.) random vectors having the same law as $(X,Y)$.
For $N\in\mathbb{N}$, set $\xi_N=(\xi^1,\ldots,\xi^N)$. We will use
approximation of $Err(f)$ by means of $\xi_N$ (as $N\to \infty$) and
a {\it prediction algorithm} (PA). This PA employs a function
$f_{PA}=f_{PA}(x,\xi_N)$ defined for $x\in\mathbb{X}$ and $\xi_N$
and taking values in $\mathbb{Y}$. More exactly, we operate with a
{\it family of functions} $f_{PA}(x,v_p)$ (with values in
$\mathbb{Y}$) defined for $x\in\mathbb{X}$ and
$v_p\in(\mathbb{X}\times\mathbb{Y})^p$ where $p\in \mathbb N$, $p\le
N$. To simplify the notation we write $f_{PA}(x,v_p)$ instead of
$f_{PA}^p(x,v_p)$. For $S\subset \{1,\ldots,N\}$ we set
$\xi_N(S)=\{\xi^j,j\in S\}$ and $\overline  S
:=\{1,\ldots,N\}\setminus S$. For  $K\in \mathbb N$ $(K>1)$,
introduce a partition of a set $\{1,\ldots,N\}$ by means of subsets
$$
S_k(N)=\{(k-1)[N/K]+1,\ldots,k[N/K]\mathbb{I}\{k<K\}+N\mathbb{I}\{k=K\}\},~k=1,\ldots,K,
$$
here $[a]$ is the integer part of a number $a\in\mathbb{R}$.
Following Bulinski (2012) we can construct an estimate of $Err(f)$
involving $\xi_N$, prediction algorithm defined by $f_{PA}$ and
$K$-cross-validation (on cross-validation we refer, e.g., to
Arlot and Celisse (2010)). Namely, set
\begin{equation}\label{Errestmult}
\widehat{Err}_K(f_{PA},\xi_N):=\sum_{i=0}^{2m-1}\sum_{i-m<
|y|\leq m} \frac{1}{K}\sum_{k=1}^K \sum_{j\in S_k(N)}
\frac{\widehat\psi(y,\xi_N({S_k(N)}))\mathbb{I}
\{A_N(y,i,k,j)\}}{\sharp
S_k(N)}
\end{equation}
where $A_N(y,i,k,j)=\{Y^j\!=\!y,|f_{PA}(X^j,\xi_N(\overline{S_k(N)}))\!-\!y|>i\}$.
Here, for each $k\in\{1,\ldots,K\}$, let
$\widehat{\psi}(y,\xi_N({S_k(N)}))$ be strongly consistent
estimates of $\psi(y)$ (as $N\to \infty$)
 for all $y\in\mathbb{Y}$, i.e.
\begin{equation*}
\widehat{\psi}(y,\xi_N({S_k(N)}))\to
\psi(y)\;\;\mbox{a.s.},\;\;y\in\mathbb{Y},\;\;N\to \infty.
\end{equation*}
In Bulinski and Rakitko (2014) the criterion was established
to guarantee the relation
\begin{equation*}
\widehat{Err}_K(f_{PA},\xi_N)\to Err(f)\;\;\mbox{a.s.},\;\;N\to \infty.
\end{equation*}
For $r\in \{1,\ldots,n\}$ set $\mathbb{X}_r=\{0,1\ldots,s\}^r$.
Then
$\mathbb{X}=\mathbb{X}_n$. We write  $\alpha =
(k_1,\ldots,k_r)$, $X_{\alpha}= (X_{k_1},\ldots,X_{k_r})$ and
$x_{\alpha}= (x_{k_1},\ldots,x_{k_r})$ where $x_i\in
\{0,\ldots,s\}$, $i=1,\ldots,n$. In many models it is natural to assume
that $Y$ depends only on some collection of factors $X_{\alpha}$.
We say that a vector $\alpha$ (and the corresponding vector $X_{\alpha}$) is {\it significant} if, for $x\in \mathbb{X}$ and $y\in\mathbb{Y}$,
one has
${\sf P}(Y=y|X=x)={\sf P}(Y=y|X_{\alpha}=x_{\alpha})$
whenever ${\sf P}(X=x)>0$. In Bulinski and Rakitko (2014) (formula (14)), for each $\beta=(m_1,\ldots,m_r)$
with $1\leq m_1 < \ldots <m_r \leq n$, the function $f^{\beta}$ was introduced
and (formula (19)) prediction algorithm $\widehat{f}^{\beta}(x,\xi_N(W_N))$ was
proposed where  $x\in\mathbb{X}$ and $\xi_N(W_N)=(\xi_{n_1},\ldots,\xi_{n_u})$,
$W_N=\{n_1,\ldots,n_u\}\subset \{1,\ldots,N\}$.
It was proved (Theorem 2 in Bulinski and Rakitko (2014)) that if $\alpha=(k_1,\ldots,k_r)$ is significant then, for any $\beta=(m_1,\ldots,m_r)$  and each $\nu >0$,
one has $\widehat{Err}_K(\widehat{f}^{\alpha}_{PA},\xi_N) \leq
\widehat{Err}_K(\widehat{f}^{\beta}_{PA},\xi_N) + \nu$ a.s. for
all $N$ large enough. Thus it is reasonable to
choose among all $\beta=(m_1,\ldots,m_r)$ such vector $\alpha$ that
$\alpha= {\sf argmin}_{\beta}\{\widehat{Err}_K(\widehat{f}^{\beta}_{PA},\xi_N)\}$
or take for further analysis (using permutation tests, see, e.g., Golland et al. (2005)) several vectors giving
the estimated prediction error close to the minimal value.
Moreover, for specified sequence $\varepsilon=(\varepsilon_N)_{N\in\mathbb{N}}$
of positive numbers, the regularized versions $\widehat{f}^{\beta}_{PA,\varepsilon}$ of
$\widehat{f}^{\beta}_{PA}$ were introduced and the CLT
was established (Theorem~3 in Bulinski and Rakitko (2014)) for these estimates.
Further extension of such CLT is obtained in Section 4 of the present paper.

\vskip 10mm

\noindent  3. SIMULATION

To illustrate our approach we consider three
examples. For each example we simulated i.i.d. random vectors
$\xi_1,\ldots,\xi_N$. Then (for each example) we evaluated the
estimate
$\widehat{Err}_K(\widehat{f}^{\beta}_{PA,\varepsilon},\xi_N)$ where
$K=10$, vector $\beta$ had appropriate dimension, and for
regularization of estimates we employed $\varepsilon_N= N^{-1/4}$,
$N\in\mathbb{N}$. After that we took all possible collections
$\beta$ of $r$ factors among $n$ and selected $10$ of them with
lowest values of estimated prediction error
$\widehat{Err}_K(\widehat{f}^{\beta}_{PA,\varepsilon},\xi_N)$. For
saving time of calculations we used $n=50$ factors. However the
results are interesting and instructive. Let the factors $X_i$,
$i=1,\ldots,n$, be i.i.d. random variables taking values $-1,0,1$
with probabilities $1/3$ and $Y$ be a binary response variable with
values $-1$ and $1$. We assume also that $r$ (the cardinality of the
collection of significant factors) is equal to 3 in Example 1 and
equals 4 in Examples 2 and 3. In Examples 1 and 2 the impact of the
``noise'' on response variable  is described by means of
multiplication of $Y$ by the random variable $(-1)^{Z_{\gamma}}$
where $Z_{\gamma}$ is the Bernoulli random variable, namely,
$\p(Z_{\gamma}=1)=\gamma$ and $\p(Z_{\gamma}=0)=1-\gamma$. We
consider $\gamma = 0.1$, that is the mean level of noise is $10\%$.
Assume that $Z_{\gamma}$ and $X=(X_1,\ldots,X_n)$ are independent.
\vskip0.5cm {\bf Example 1}. Let $r=3$ and $Y=Y^{0}\cdot
(-1)^{Z_{\gamma}}$ where
$$
Y^{0}=
\begin{cases}
\;\;1,~~~&X_2=1,~X_3\ge0,
\vspace{-0.3cm}
\\
\;\;1,&X_2=-1,~X_3+X_5\ge 1,
\vspace{-0.3cm}\\
-1,~~~&\text{otherwise}.
\end{cases}
$$
Here $X_2, X_3, X_5$ are the  factors determining $Y$.
\vskip0.5cm

{\bf Example 2}. Take $r=4$ and set $Y=Y^{0}\cdot (-1)^{Z_{\gamma}}$ where
$$
Y^0=
\begin{cases}
\;\;1,~~~&X_2=1,
\vspace{-0.3cm}
\\
\;\;1,&X_3+X_5+X_8\ge 2,
\vspace{-0.3cm}
\\
-1,~~~&\text{otherwise}.
\end{cases}
$$
The  factors determining $Y$ are  $X_2, X_3, X_5, X_8$.
\vskip0.5cm
In the following example we consider nonlinear constrains.

{\bf Example 3}. Let $r=4$. Set
$$
Y=
\begin{cases}
\;\;1,~~~&3^{X_1+X_2+X_4}\sin(X_3Z^{\ln(X_3-2X_4+7)})>1,
\vspace{-0.3cm}
\\
-1,~~~&\text{otherwise,}
\end{cases}
$$
assuming the random variable
$Z$ be uniformly distributed on $[0,1]$. Let $Z$ and  $X$ be
independent. Here $X_1, X_2, X_3, X_4$ are the factors determining
$Y$.

\vskip0.5cm Collections of various factors and corresponding  values
of $\widehat{Err}_K(\widehat{f}^{\beta}_{PA,\varepsilon},\xi_N)$
obtained for $N=500$ are presented in Tables 1, 2 and 3. Namely,
$EPE_i$  stands for $\widehat{Err}_K$ found in the framework of
Example $i$ where $i=1,2,3$. Columns $n_1, n_2, n_3$ (and $n_1, n_2,
n_3, n_4$) in the tables indicate the choice of factors $X_{n_1},
X_{n_2}, X_{n_3}$ (and $X_{n_1}, X_{n_2}, X_{n_3}, X_{n_4}$),
respectively. The same information is provided in Tables 4, 5 and 6
where one has $N=1000$.

It is worth to emphasize that in all considered examples for large
($N=1000$) and rather modest ($N=500$) samples our method permits to
identify correctly the collections of
significant factors (corresponding to the minimum of prediction
error estimates). Moreover, these tables show that the estimated
prediction error for significant collections of factors has visible
advantage w.r.t. other collections.

\vskip1cm

\begin{multicols}{3}
\begin{table}[ht]
\centering
\begin{tabular}{ccc|c}
  \hline
$n_1$ & $n_2$ & $n_3$ & $EPE_1$ \\
  \hline
   {\bf 2} & {\bf 3} & {\bf 5} & {\bf 0.6336 }\\
   2 &  3 & 32 & 0.8020 \\
   2 &  3 & 48 & 0.8100 \\
   2 &  3 & 28 & 0.8260 \\
   2 &  3 &  4 & 0.8515 \\
   2 &  3 & 31 & 0.8527 \\
   2 &  3 & 22 & 0.8528 \\
   2 &  3 & 34 & 0.8551 \\
   2 &  3 & 50 & 0.8649 \\
   2 &  3 & 23 & 0.8652 \\
   \hline
\end{tabular}
\caption{$r=3$, N=500}
\label{tab1}
\end{table}

\begin{table}[ht]
\centering
\begin{tabular}{cccc|c}
  \hline
$n_1$ & $n_2$ & $n_3$ & $n_4$ & $EPE_2$ \\
  \hline
 {\bf 2} & {\bf 3} &  {\bf 5} &  {\bf 8} & {\bf 0.3997} \\
   2 &  3 &  5 & 24 & 0.5901 \\
   2 &  3 &  5 & 46 & 0.5911 \\
   2 &  3 &  5 & 32 & 0.5961 \\
   2 &  3 &  5 & 31 & 0.6014 \\
   2 &  3 &  5 & 10 & 0.6059 \\
   2 &  3 &  5 & 14 & 0.6224 \\
   2 &  3 &  5 & 42 & 0.6250 \\
   2 &  3 &  5 & 29 & 0.6251 \\
   2 &  3 &  5 & 22 & 0.6267 \\
   \hline
\end{tabular}
\caption{$r=4$, N=500}
\label{tab2}
\end{table}

\begin{table}[ht]
\centering
\begin{tabular}{cccc|c}
  \hline
$n_1$ & $n_2$ & $n_3$ & $n_4$ & $EPE_3$ \\
  \hline
 {\bf 1} &  {\bf 2} & {\bf 3} &  {\bf 4} & {\bf 0.0939} \\
   1 &  3 & 20 & 42 & 0.2956 \\
   2 &  3 &  5 & 29 & 0.3211 \\
   1 &  3 &  4 & 39 & 0.3228 \\
   1 &  2 &  3 &  8 & 0.3322 \\
   1 &  3 & 24 & 42 & 0.3355 \\
   1 &  2 &  3 &  5 & 0.3395 \\
   1 &  2 &  3 & 20 & 0.3431 \\
   1 &  2 &  3 & 40 & 0.3487 \\
   1 &  2 &  3 & 27 & 0.3558 \\
   \hline
\end{tabular}
\caption{$r=4$, N=500}
\label{tab3}
\end{table}
\end{multicols}

\newpage

\begin{multicols}{3}
\begin{table}[ht]
\centering
\begin{tabular}{ccc|c}
  \hline
$n_1$ & $n_2$ & $n_3$ & $EPE_1$ \\
  \hline
 {\bf 2} &  {\bf 3} & {\bf  5} & {\bf 0.5675} \\
  2 &  3 & 32 & 0.7981 \\
   2 &  3 & 47 & 0.8096 \\
   2 &  3 & 34 & 0.8126 \\
   2 &  3 &  4 & 0.8127 \\
   2 &  3 & 44 & 0.8334 \\
   2 &  3 & 48 & 0.8369 \\
   2 &  3 & 22 & 0.8401 \\
   2 &  3 & 23 & 0.8441 \\
   2 &  3 & 31 & 0.8442 \\
   \hline
\end{tabular}
\caption{$r=3$, N=1000}
\label{tab4}
\end{table}

\begin{table}[ht]
\centering
\begin{tabular}{cccc|c}
  \hline
$n_1$ & $n_2$ & $n_3$ & $n_4$ & $EPE_2$ \\
  \hline
 {\bf 2} & {\bf 3} & {\bf 5} & {\bf 8} & {\bf 0.4768} \\
   2 &  3 &  5 & 42 & 0.6936 \\
   2 &  5 &  8 & 11 & 0.6970 \\
   2 &  5 &  8 & 26 & 0.6974 \\
   2 &  3 &  5 &  6 & 0.6981 \\
   2 &  5 &  8 & 12 & 0.7035 \\
   2 &  5 &  8 & 50 & 0.7039 \\
   2 &  3 &  5 & 32 & 0.7045 \\
   2 &  3 &  5 & 27 & 0.7060 \\
   2 &  3 &  5 & 46 & 0.7063 \\
   \hline
\end{tabular}
\caption{$r=4$, N=1000}
\label{tab5}
\end{table}

\begin{table}[ht]
\centering
\begin{tabular}{cccc|c}
  \hline
$n_1$ & $n_2$ & $n_3$ & $n_4$ & $EPE_3$ \\
  \hline
{\bf 1} & {\bf 2} & {\bf 3} & {\bf 4} & {\bf 0.2278} \\
   2 &  3 &  4 & 32 & 0.3355 \\
   1 &  2 &  3 &  6 & 0.4352 \\
   1 &  2 &  3 & 46 & 0.4663 \\
   1 &  3 &  4 & 15 & 0.4694 \\
   1 &  2 &  3 & 27 & 0.4697 \\
   1 &  2 &  3 & 50 & 0.4704 \\
   2 &  3 &  4 & 18 & 0.4812 \\
   2 &  3 &  4 & 44 & 0.4856 \\
   1 &  3 &  4 & 40 & 0.4862 \\
   \hline
\end{tabular}
\caption{$r=4$, N=1000}
\label{tab6}
\end{table}

\end{multicols}

\vskip1cm

However, if $N$ is not large enough the proposed stochastic approach
can lead to the choice of a collection of factors which is not (the
most) significant. For instance, if $N=500$ then the right
identifications of significant factors
have occurred  in $99\%$, $97\%$, $69\%$ of respective simulations
for
Examples 1, 2 and 3
 (averaging is over 100 performance procedures). In Example 3 this
frequency of right identification
increases till
$93\%$ when
$N=1000$.

Figures 1, 2 and 3 demonstrate for each example the character of
stabilization of $\widehat{Err}_K$ fluctuations as $N$ grows. This
stabilization of estimates can be explained  not only by their
strong consistency but also on account of
their asymptotic normality.
In this regard we concentrate further on the new conditions which
guarantee the CLT validity for proposed prediction error estimates.

\begin{figure}
    \centering

\includegraphics[width=\textwidth,height=\textheight,keepaspectratio]{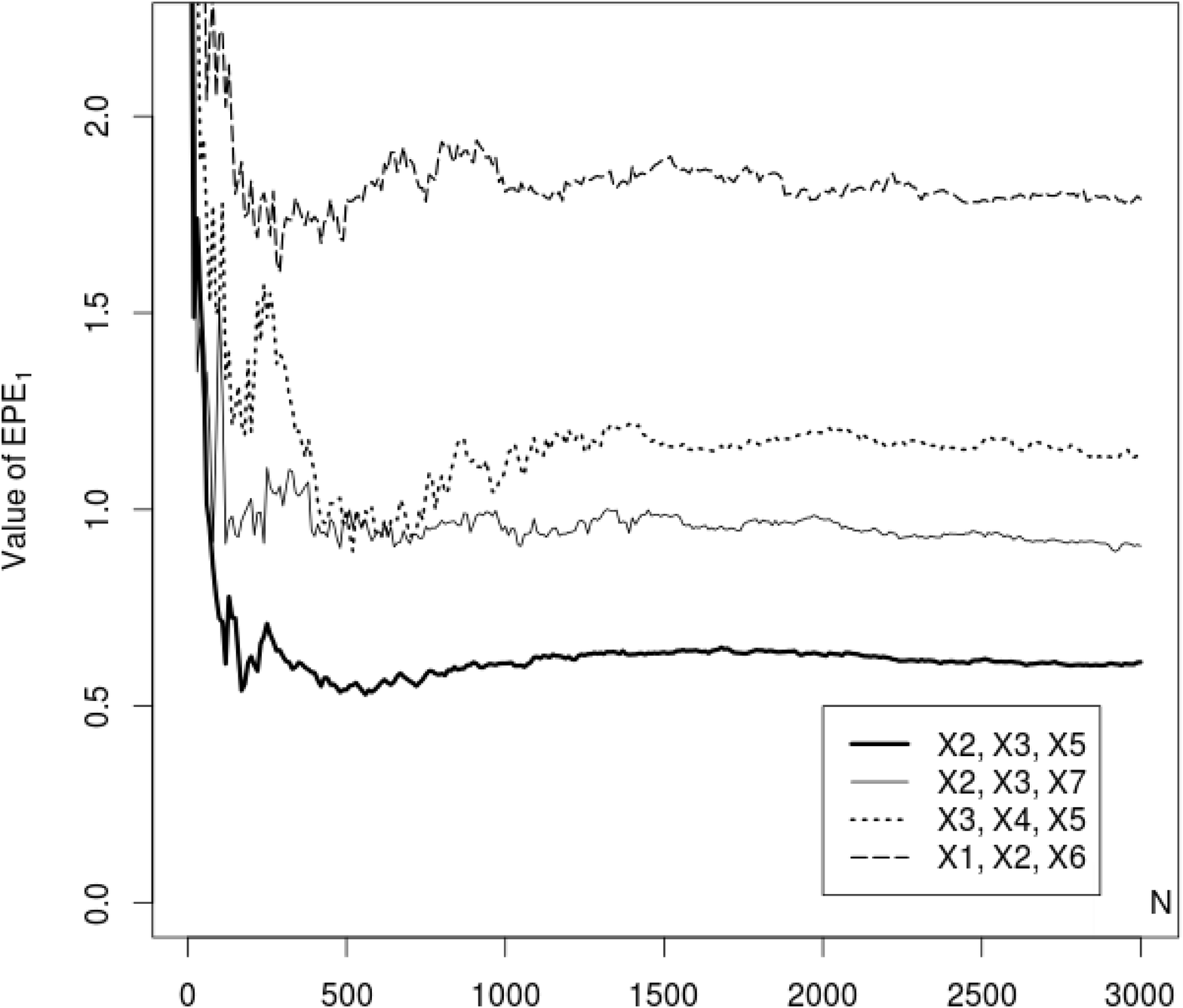}
\caption{Simulations corresponding to Example 1.}

\includegraphics[width=\textwidth,height=\textheight,keepaspectratio]{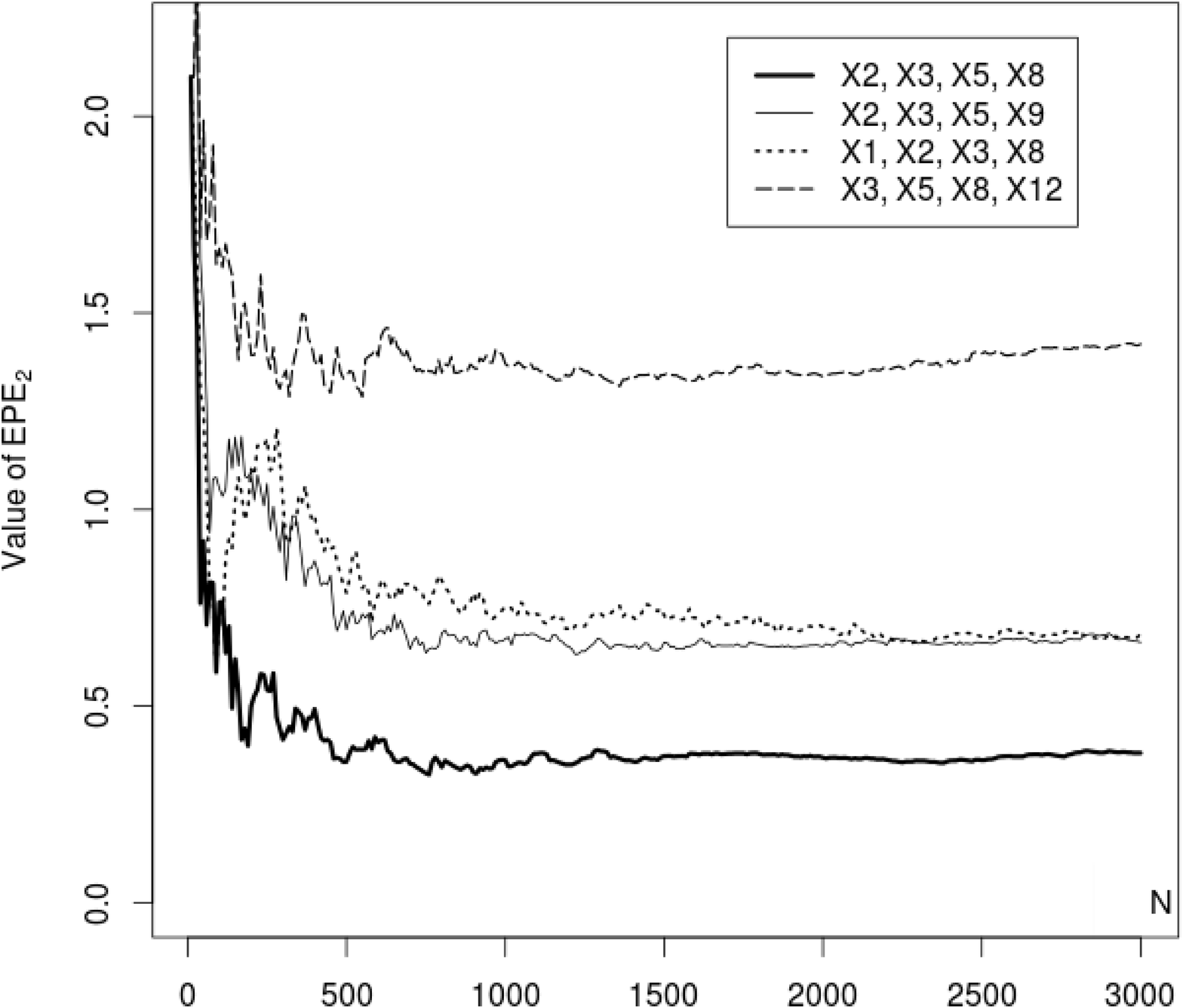}
\caption{Simulations corresponding to Example 2.}

\includegraphics[width=\textwidth,height=\textheight,keepaspectratio]{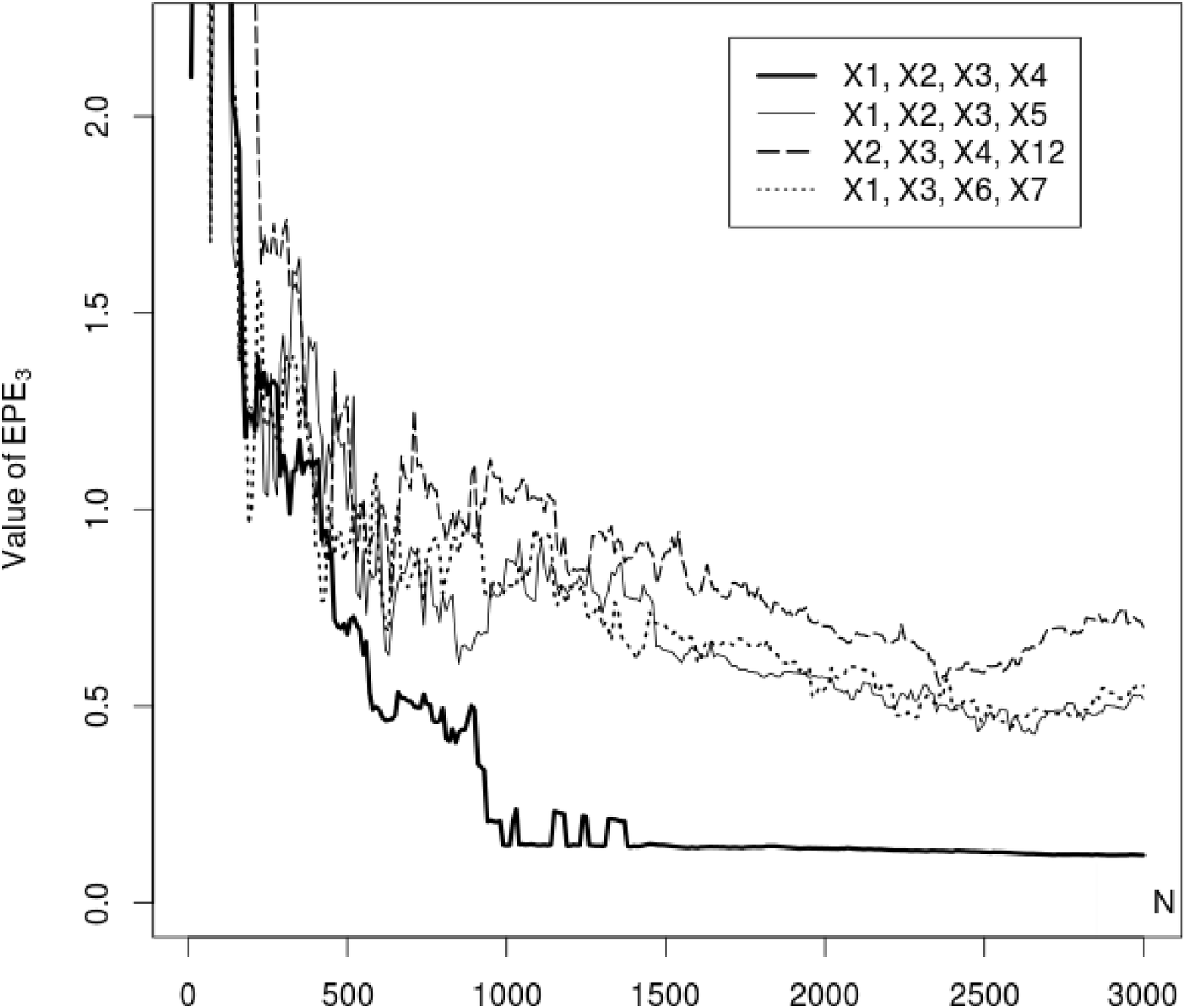}
\caption{Simulations corresponding to Example 3.}

\end{figure}

\vskip6cm

\noindent 4. NEW VERSION OF THE CENTRAL LIMIT THEOREM

We proved in Bulinski and Rakitko (2014) that asymptotic distribution of
random variables $\sqrt{N}(\widehat{Err}_K(f_{PA},\xi_N)-Err(f))$ coincides with the limit law of
\begin{equation}\label{eq_6}
\sqrt{N}(\widehat{T}_N(f)-Err(f))
=\dfrac{\sqrt{N}}{K}\sum_{k=1}^K\sum_{i=0}^{2m-1}\sum_{i-m<|y|\le m}
\dfrac{1}{\sharp S_k(N)}\sum_{j\in S_k(N)}
h_N(y,i,k,j),
\end{equation}
as $N\to\infty$, where
\begin{equation*}
h_N(y,i,k,j)= \widehat\psi(y,{S_k(N)})\ind\{Y^j=y,|f(X^j)-y|>i\}-\psi(y)\p(Y=y,|f(X)-y|>i)
\end{equation*}
and
$\widehat\psi(y,{S_k(N)}):=\widehat\psi(y,\xi_N(S_k(N)))$.

Evidently the summands here are not independent in view of the
presence of $\widehat{\psi}(\cdot,{S_k(N)})$. To prove the CLT for
random variables appearing in \eqref{eq_6} we used in
Bulinski and Rakitko (2014) the hypothesis of asymptotic normality of the
vector consisting of two subvectors, one of them being
$\sqrt{N}(\widehat{\psi}(\cdot,{S_k(N)})-\psi(\cdot))$.
Now we employ another approach assuming
symmetry of the estimates
$\widehat{\psi}(\cdot,{S_k(N)})$ of a penalty function.
Recall the following

\begin{df}
A collection of random variables $(X_1,\ldots,X_n)$, $n\in \mathbb N$, is called {\it exchangeable} if, for any permutation $\sigma\in S(n)$ of the set $\{1,\ldots,n\}$,
one has
\begin{equation*}
Law(X_1,\ldots,X_n) = Law(X_{\sigma(1)},\ldots, X_{\sigma(n)}).
\end{equation*}
\end{df}

Take $K\in \mathbb{N}$
and suppose  that $N/K=q$ where $q\in\mathbb{N}$.
Thus $\sharp S_k(N)=q$ for each $k=1,\ldots,K$.
Consider the sequence of $K\times q$ matrices $(C^{(N)})_{N\in\mathbb{N}}$ with entries
\begin{equation}\label{xi_N}
\xi^{(N)}_{k,j}:=\sum_{i=0}^{2m-1}\sum_{i-m<|y|\le m}\widehat{\psi}(y,S_k(N))\cdot\ind\{Y^{j+(k-1)q}=y, |f(X^{j+(k-1)q})-y|>i\}
\end{equation}
where $k=1,\ldots,K$ and $j=1,\ldots,q$.
Introduce
\begin{equation}\label{trarr}
X_{N,j}:=\dfrac{1}{\sqrt{K}}\sum_{k=1}^K\xi_{k,j}^{(N)},~~~j=1,\ldots,q.
\end{equation}
Then
\begin{equation}\label{sum_main}
\sqrt{N}(\widehat{T}_N(f)-Err(f))=\frac{1}{\sqrt{q}}\sum_{j=1}^q (X_{N,j} - \sqrt{K}Err(f)).
\end{equation}

We take the functions $\widehat{\psi}(y,\cdot)$ which are symmetric for each $y\in \mathbb Y$. Then any row and any column of $C^{(N)}$ contain
exchangeable random variables (row-column exchangeability).
Clearly, the triangular array $\{X_{N,j},~1\le j\le q, N\in\mathbb{N}\}$
is row-wise exchangeable.

We will establish
the CLT for sums appearing in \eqref{sum_main}.
In Berti et al. (2004) one can find several results which guarantee the CLT
validity when the summands  $\{X_i\}_{i=1}^n$ are (in appropriate manner) conditionally identically distributed.
Namely,
\begin{equation}\label{clt_type}
\dfrac{1}{\sqrt{n}}\big(f(X_1)+\ldots+f(X_n)-L_n\big)\xrightarrow{law}Z_{0,\sigma^2}\sim\mathcal{N}(0,\sigma^2)
\end{equation}
where $f$ is a measurable function such that $\E|f(X_1)|<\infty$
and $L_n=L_n(X_1, \ldots, X_n)$. In the mentioned paper the authors
applied the martingale techniques. Such approach was developed for
exchangeable variables in Weber (1980).
We will prove the CLT in the form \eqref{clt_type} with $f(x)=x$ for row-wise exchangeable arrays  by means of other tools.
We will employ the recent result of R{\"o}llin (2013).
Let ${\sf Y}=(Y_1,\ldots,Y_m)$ be a collection of exchangeable random variables such that
\begin{equation}\label{cond1}
\E Y_{1} = 0,~~~\E|Y_{1}|^3<\infty.
\end{equation}
Consider $\Sigma=(\sigma_{i,j})_{1\le i,j\le m}$ with
$\sigma_{i,j}=\E(Y_iY_j)$, i.e. the covariance matrix of $\sf Y$.
Set $\sigma_{i,i}=\sigma^2$. Suppose that $Y_1+ \ldots + Y_m =C_m$
a.s. where $C_m$ is a constant. Then w.l.g. we can assume that
\begin{equation}\label{cond2}
 \sum_{i=1}^mY_i=0~~~{a.s.}
\end{equation}
For a function $h:\mathbb R^d\rightarrow \mathbb R$
and $k\in\mathbb{N}$ set
$$
C_h^{(k)}:=\max_{i_1,\ldots,i_d\geq 0, \sum_{j=1}^d i_j=k}\left\| \frac{\partial^kh}{\partial x_1^{i_1}\ldots \partial x_d^{i_d}}\right\|_{\infty}.
$$

\begin{Th}[R{\"o}llin (2013)]\label{Rollin}
Let $\sf Y$ be a vector consisting of exchangeable random variables
and having a covariance matrix $\Sigma$.
Assume that conditions \eqref{cond1} and \eqref{cond2} are satisfied. Then
\begin{equation}\label{stein}
|\E h({\sf Y})-\E h({\sf Z})|\le C_h^{(2)}\Big[{\sf var}\Big(\sum_{i=1}^mY_i^2\Big)\Big]^{\frac{1}{2}}+16mC_h^{(3)}\E |Y_1|^3
\end{equation}
where ${\sf Z} \sim \mathcal{N}(0,\Sigma)$.
\end{Th}

For an array $\{X_{n,i}, 1\le i \le k_n, n\in\mathbb{N}\}$
we will use the following notation
\begin{equation}\label{mu}
\widehat \mu_{k_n}:=\dfrac{1}{k_n}\sum_{i=1}^{k_n}X_{n,i},\;\;\;\;
\widehat \sigma^2_{k_n}:=\dfrac{1}{k_n}\sum_{i=1}^{k_n}(X_{n,i}-\widehat\mu_{k_n})^2.
\end{equation}

We apply \eqref{stein} to prove the following result.

\begin{Lem}\label{Lem1}
Let $\{X_{n,i}, 1\le i \le k_n, n\in\mathbb{N}\}$
be a row-wise exchangeable array where positive integers $k_n\to \infty$ as $n\to \infty$. Suppose that
\vskip0.1cm
\;\;
$
1^{\circ}\!.\;\; \sup_{n\in\mathbb{N}}\E X_{n,1}^4<\infty,
$
\vskip0.1cm
\;\;
$
2^{\circ}\!.\;\;\E X_{n,1}^2-\E X_{n,1}X_{n,2}\to\sigma^2>0,\;n\to\infty,
$
\vskip0.1cm
\;\;
$
3^{\circ}\!.\;\;{\sf cov}(X_{n,1}^2,X_{n,2}^2)+{\sf cov}(X_{n,1}X_{n,2},X_{n,3}X_{n,4})-2\,{\sf cov}(X_{n,1}^2,X_{n,2}X_{n,3})\to0,\;n\to \infty.
$
\vskip0.1cm
\noindent
Then, for any sequence $(m_n)_{n\in\mathbb{N}}$ of positive integers
such that $m_n\to\infty$ and $m_n/k_n\to \alpha <1$  as $n\to \infty$, the following relation holds
\begin{equation*}
\dfrac{1}{\sqrt{m_n}}\sum_{i=1}^{m_n}(X_{n,i}-\widehat\mu_{k_n})
\xrightarrow{law}Z_{0,(1-\alpha)\sigma^2}\sim\mathcal{N}(0,(1-\alpha)\sigma^2),\;\;n\to\infty.
\end{equation*}
\end{Lem}
{\it Proof}.
First of all, for each $n\in\mathbb{N}$, we introduce the auxiliary random variables
$$
Y_{n,i}:=X_{n,i}-\widehat\mu_{k_n},~~~i=1,\ldots,k_n.
$$
The collection $\{Y_{n,1},\ldots, Y_{n,k_n}\}$ is exchangeable as
$\{X_{n,1},\ldots,X_{n,k_n}\}$ has this property. Obviously $\sum_{i=1}^{k_n}Y_{n,i}=0$ a.s.
for any $n\in\mathbb{N}$.
Moreover, $\E Y_{n,1}=0,$ for any $n\in\mathbb{N}$.
One can verify that
$$
\E Y_{n,1}^2
=\left(1-\frac{1}{k_n}\right)(\E X_{n,1}^2-\E X_{n,1}X_{n,2}),\;\;\;
\E Y_{n,1}Y_{n,2}
=-\,\frac{1}{k_n}(\E X_{n,1}^2-\E X_{n,1}X_{n,2}).
$$
For each $n\in \mathbb{N}$, take a vector ${\bf Z}=(Z_{n,1},\ldots,Z_{n,m_n})$ independent of $(X_{n,1},\ldots,X_{n,k_n})$
and such that ${\bf Z}\sim \mathcal{N}(0,\Sigma)$. Here $\Sigma$ is a covariance matrix of ${\bf Y}=(Y_{n,1},\ldots, Y_{n,m_n})$. Thus ${\sf cov}(Z_{n,i},Z_{n,j})={\sf cov}(Y_{n,i},Y_{n,j})$, $1\le i,j\le m_n$. Clearly,
$$\dfrac{1}{\sqrt{m_n}}\sum_{i=1}^{m_n}(X_{n,i}-\widehat\mu_{k_n})=
\frac{1}{\sqrt{m_n}}\sum_{i=1}^{m_n}Y_{n,i}=:
S_{{\bf Y}, m_n}.$$
Set $S_{{\bf Z},m_n}:=\frac{1}{\sqrt{m_n}}\sum_{i=1}^{m_n}Z_{n,i}.$
In view of $2^{\circ}$ condition $m_n/k_n\to \alpha$ ($n\to \infty$) yields
$${\sf var} S_{{\bf Z},m_n}= \E Y_{n,1}^2 + (m_n-1)\E Y_{n,1}Y_{n,2}=\left(1-\dfrac{m_n}{k_n}\right)(\E X_{n,1}^2-\E X_{n,1}X_{n,2})\to (1-\alpha)\sigma^2.$$
Consequently, $ S_{{\bf Z},m_n}\stackrel{law}\longrightarrow
\mathcal{N}(0,(1-\alpha)\sigma^2),\;\;n\to \infty. $ Now we show
that $S_{{\bf Y},m_n}$ and $S_{{\bf Z},m_n}$ have the same limit
distribution. Due to Theorem 7.1 Billingsley (1968) it is sufficient
to verify that
\begin{equation}\label{conv1}
\E f(S_{{\bf Y},m_n}) - \E f(S_{{\bf Z},m_n})\to 0,~~~ n\to\infty,
\end{equation}
for any three times continuously differentiable function $f:\mathbb R \to \mathbb R$ such that
$$
c_{f}^{(j)}:= \left\|\frac{d^jf}{dx^j}\right\|_{\infty} <\infty,\;\;j=1,2,3.
$$
For any fixed $n\in\mathbb{N}$, apply Theorem \ref{Rollin} with $m=m_n$, $Y_i=\frac{1}{\sqrt{m_n}}Y_{n,i}$, $i=1,\ldots,m_n$, and
$$
h(x_1,\ldots,x_{m_n}):=f(x_1+\ldots+x_{m_n}).
$$
Then we can write
$$
\big|\E f(S_{{\bf Y},m_n}) - \E f(S_{{\bf Z},m_n})\big|=\big|\E h({\bf Y}) - \E h({\bf Z})\big|
$$
$$
 \le C_f^{(2)}m_n^{-1}\Big[{\sf var}\Big(\sum_{i=1}^{m_n}Y_{n,i}^2\Big)\Big]^{\frac{1}{2}}+
 16 C_f^{(3)} m_n^{-1/2}\E |Y_{n,1}|^3.
 $$
Note that
$$
{\sf var}\Big(\sum_{i=1}^{m_n}Y_{n,i}^2\Big)=m_n\E Y_{n,1}^4+m_n(m_n-1)\E Y_{n,1}^2Y_{n,2}^2-m_n^2\big(\E Y_{n,1}^2\big)^2
$$
$$
=m_n\big(\E Y_{n,1}^4-(\E Y_{n,1}^2)^2\big)+m_n(m_n-1){\sf cov}(Y_{n,1}^2,Y_{n,2}^2).
$$
We claim that
$${\sf cov}(Y_{n,1}^2,Y_{n,2}^2)-\Big[{\sf cov}(X_{n,1}^2,X_{n,2}^2)+{\sf cov}(X_{n,1}X_{n,2},X_{n,3}X_{n,4})-2{\sf cov}(X_{n,1}^2,X_{n,2}X_{n,3})\Big]\to0$$
as $n\to\infty$.
Indeed, set $S_n=\frac{1}{k_n}\sum_{i=1}^{k_n}X_{n,i}$. Using exchangeability property
of $(X_{n,1},\ldots,X_{n,k_n})$ and taking into account that
covariance function is bilinear we obtain
$$
\Co(Y_{n,1}^2,Y_{n,2}^2)=\E Y_{n,1}^2Y_{n,2}^2-(\E Y_{n,1}^2)^2
$$
$$
=\Co(X_{n,1}^2,X_{n,2}^2)+2\,\Co(X_{n,1}^2,S_n^2)-4\,\Co(X_{n,1}^2,X_{n,2}S_n)
$$
$$
-4\,\Co(X_{n,1}S_n,S_n^2)+4\,\Co(X_{n,1}S_n,X_{n,2}S_n)+\Co(S_n^2,S_n^2).
$$
For $n\to \infty$, by virtue of $1^{\circ}$ we get
$$
\Co(X_{n,1}^2,S_n^2)=\Co (X_{n,1}^2,X_{n,2}X_{n,3})+O(k_n^{-1}),
$$
$$
\Co(X_{n,1}^2,X_{n,2}S_n)=\Co(X_{n,1}^2,X_{n,2}X_{n,3})+O(k_n^{-1}),
$$
$$
\Co(X_{n,1}S_n,S_n^2)=\Co (X_{n,1}X_{n,2},X_{n,3}X_{n,4})+O(k_n^{-1}),
$$
$$
\Co (X_{n,1}S_n,X_{n,2}S_n)=\Co (X_{n,1}X_{n,2},X_{n,3}X_{n,4})+O(k_n^{-1}).
$$
Therefore, condition $3^{\circ}$ implies that ${\sf cov}(Y_{n,1}^2,Y_{n,2}^2)\to0$ as $n\to \infty$.
Thus relation \eqref{conv1} holds and the proof is complete. $\square$

\begin{Rem}
Assume that
$$\sup_{n\in\mathbb N}\E \big((X_{n,1}-\widehat \mu_{k_n})\slash\widehat\sigma_{k_n}\big)^4<\infty.
$$
Then, for a sequence $(m_n)_{n\in\mathbb{N}}$ appearing in Lemma \ref{Lem1},
one can prove the following version of the CLT
$$
\dfrac{1}{\sqrt{m_n}}\sum_{i=1}^{m_n}\Big(\dfrac{X_{n,i}-\widehat\mu_{k_n}}{\widehat
\sigma_{k_n}}\Big)\stackrel{law}\longrightarrow Z_{0,1-\alpha}\sim \mathcal{N}(0,1-\alpha),~~~n\to\infty.
$$
\end{Rem}

\begin{Rem}
In Chernoff and Teicher (1958) the result similar to Lemma 1 was established but
the important case $\alpha=0$ (which we consider further) was not
comprised. One can also employ the martingale approach of
Weber  (1980) to obtain the result of Lemma 1. However Rollin's
Theorem \ref{Rollin}  permits us to estimate the convergence rate to
the limit Gaussian law.
Moreover, we can prove that
under certain conditions the asymptotic behavior of the specified
partial sums
is described by the mixture of the normal laws.
\end{Rem}

Now we  consider the triangular array  $\{X_{N,i}, 1\le i \le q, N\in\mathbb{N}\}$
with elements defined by \eqref{trarr}.
Thus we take $k_n=q$ in Lemma 1 and write $N$ instead of $n$.

\begin{Lem}\label{Lem2}
Suppose that, for each $N\in\mathbb N$, any $y\in\mathbb Y$
and all $k=1,\ldots,K$,
\begin{equation}\label{eta}
\sup_{y\in\mathbb Y,\;N\in\mathbb N,\;k\in\{1,\ldots,K\}}\E \left(\widehat{\psi}(y,S_k(N))\right)^4<\infty.
\end{equation}
Let $(m_N)_{N\in\mathbb{N}}$
be a sequence of positive integers such that
$m_N \leq q$, $m_N\to \infty$ and\break $m_N/N \to \alpha <1$ as $N\to \infty$.
Then
\begin{equation*}
\dfrac{1}{\sqrt{m_N}}\sum_{i=1}^{m_N}(X_{N,i}-\widehat\mu_{N})\xrightarrow{law}
Z_{0,(1-\alpha)\sigma^2}\sim\mathcal{N}(0,(1-\alpha)\sigma^2)
\end{equation*}
where $\mu_N$ is introduced in \eqref{mu} $($with $k_n=q$ and $n$ replaced by $N)$ and
\begin{equation}\label{sigma^2}
\sigma^2\!=\!\E\left[\sum_{i=0}^{2m-1}\sum_{i-m<|y|\le m}\psi(y)(\ind\{Y\!=\!y, |f(X)-y|\!>\!i\}-\p(Y\!=\!y,|f(X)-y|\!>\!i))\right]^2\!\!.
\end{equation}
\end{Lem}

{\it Proof.} We show that conditions of Lemma  \ref{Lem1} are met.
 $1^{\circ}$
follows by virtue of \eqref{eq_6}, \eqref{trarr} and \eqref{eta} as
indicator function takes values in the set $\{0,1\}$. Now we turn to
$2^{\circ}$. The exchangeability of the columns of the array
$\{\xi^{(N)}_{k,j}\}$ implies that
$$
\E X_{N,1}X_{N,2}=\dfrac{1}{K}\E \left(\sum_{k=1}^K\xi^{(N)}_{k,1}\right)\left(\sum_{k=1}^K \xi^{(N)}_{k,2}\right)=\E\xi_{1,1}^{(N)}\xi_{1,2}^{(N)}+ (K-1)\E\xi^{(N)}_{1,1}\xi^{(N)}_{2,2}.
$$
The Lebesgue theorem on majorized convergence yields
that the limit behavior of $\E\xi_{1,1}^{(N)}\xi_{1,2}^{(N)}$ as $N\to \infty$
will be the same as for $\E\zeta_{1,1}^{(N)}\zeta_{1,2}^{(N)}$ where
\begin{equation*}
\zeta_{k,j}^{(N)}:=\sum_{i=0}^{2m-1}\sum_{i-m<|y|\le m}\psi(y)\,\ind\{Y^{j+(k-1)q}=y, |f(X^{j+(k-1)q})-y|>i\}.
\end{equation*}
Random vectors $(X^1,Y^1), (X^2,Y^2),\ldots$ are independent. Therefore,
$\E\zeta_{1,1}^{(N)}\zeta_{1,2}^{(N)}=\E\zeta_{1,1}^{(N)}
\E\zeta_{1,2}^{(n)}$ and in view of \eqref{Err} we get
$$
\lim_{N\to\infty}\E\xi^{(N)}_{1,1}\xi^{(N)}_{1,2}=\lim_{N\to\infty}
\big(\E\zeta_{1,1}^{(N)}\big)^2=\big(Err(f)\big)^2.
$$
In a similar way we come to the relation
$$
\lim_{N\to\infty}\E\xi^{(N)}_{1,1}\xi^{(N)}_{2,2}=\lim_{N\to\infty}
\big(\E\zeta_{1,1}^{(N)}\big)^2=\big(Err(f)\big)^2.
$$
Thus $\E X_{N,1} X_{N,2}\to K\big(Err(f)\big)^2$ as $N\to\infty$.
Applying the Lebesgue theorem once again we conclude that
$$
\lim_{N\to\infty}\E(X_{N,j})^2=\lim_{N\to\infty}\E(Z_{N,j})^2=\lim_{N\to\infty} \Big[\E(\zeta^{(N)}_{1,1})^2+(K-1)\E\zeta^{(N)}_{1,1}\zeta^{(N)}_{1,2} \Big],$$
where
\begin{equation*}
Z_{N,j}:=\dfrac{1}{\sqrt{K}}\sum_{k=1}^K\zeta_{k,j}^{(N)},~~~j=1,\ldots,q.
\end{equation*}
Taking into account that $\E\zeta^{(N)}_{1,1}\zeta^{(N)}_{1,2}=\big(Err(f)\big)^2$
(for each $N \geq 2K$) we get
$$
\sigma^2=\E\Big[\sum_{i=0}^{2m-1}\sum_{i-m<|y|\le m}\psi(y)\,\ind\{Y=y, |f(X)-y|>i\}\Big]^2-\big(Err(f)\big)^2
$$
$$
=\E\left[\sum_{i=0}^{2m-1}\sum_{i-m<|y|\le m}\psi(y)\big(\ind\{Y=y, |f(X)-y|>i\}-\p(Y=y,|f(X)-y|>i)\big)\right]^2\!.
$$
To complete the proof we verify condition $3^{\circ}$.
Due to the Lebesgue theorem
$$
\lim_{N\to \infty}{\sf Cov}(X_{N,1}^2,X_{N,2}^2)=\lim_{k\to\infty}{\sf Cov}(Z_{k,1}^2,Z_{k,2}^2)=0
$$
as $Z_{k,1}$ and $Z_{k,2}$ are independent. Quite similar arguments
justify the following relations
${\sf Cov}(X_{N,1}X_{N,2},X_{N,3}X_{N,4})\to0$ and ${\sf Cov}(X_{N,1}^2,X_{N,2}X_{N,3})\to0$ as $N\to\infty$. $\square$

Let us discuss the established result. Instead of the initial task
to study asymptotic behavior of $\sqrt{N}(\widehat{T}_N(f)-Err(f))$
we are able to specify the limit law for difference of two estimates of $Err(f)$.
Namely, set
$$
\widehat{L}_{m_N}=\dfrac{1}{m_N}\sum_{j=1}^{m_N}\dfrac{1}{\sqrt{K}}\sum_{k=1}^K\sum_{i=0}^{2m-1}\sum_{i-m<|y|\le m}\Big[\widehat{\psi}(y,S_k(N))\cdot\ind\{Y^{j+(k-1)q}=y, |f(X^{j+(k-1)q})-y|>i\}\Big]
$$
and introduce $\widehat {L}_q$ by the same formula with $q$ instead of $m_n$. Then Lemma \ref{Lem2} affirms that
$\sqrt{m_N}(\widehat{L}_{m_N}-\widehat{L}_q)\stackrel{law}\longrightarrow
Z_{0,(1-\alpha)\sigma^2}\sim\mathcal{N}(0,(1-\alpha)\sigma^2)$ as $N\to\infty$.
Therefore, if we provide conditions to
guarantee that
$\sqrt{m_N}\big(\widehat{L}_{q}-Err(f)\big)\xrightarrow{{\sf P}}0$
then we can construct the approximate confidence intervals for $Err(f)$.
We demonstrate that this is possible for regularized
statistics introduced in Bulinski and Rakitko (2014)
to identify the significant collections of factors.

For a sequence of random variables $(\eta_N)_{N\in\mathbb{N}}$ we
write $\eta_N=O_{\sf P}(1)$ if, for any $\gamma >0$, there exists
$M(\gamma)>0$ such that ${\sf P}(|\eta_N|\geq M(\gamma)) \leq
\gamma$ for all $N$ large enough. Let $(m_N)_{N\in\mathbb{N}}$ be a
sequence of positive integers such that
$m_N\leq q$ for $q=[N/K]$ and
\begin{equation*}
m_N\to \infty,\;\;m_N/N\to 0,\;\;\mbox{as}\;\;N\to \infty.
\end{equation*}

\begin{Th}\label{Th2} Let $(m_N)_{N\in\mathbb{N}}$ be a sequence introduced above. Assume that
$\varepsilon=(\varepsilon_N)_{N\in\mathbb{N}}$ is a sequence of
positive numbers such that $\varepsilon_N\to 0$ and
$m_N^{1/2}\varepsilon_N \to \infty$ as $N\to \infty$. Take any
vector $\beta=(k_1,\ldots k_r)$ with $1\!\leq\! k_1\!<\!\ldots \!<\!
k_r\!\leq\! n$, the corresponding function $f=f^\beta$ and the
prediction algorithm defined by
$f_{PA}=\widehat{f}_{PA,\varepsilon}^{\beta}$. Let, for any
$y\in\mathbb Y$  and $k\in\{1,\ldots,K\}$, the estimate
$\widehat{\psi}(y,S_k(N))$ be strongly consistent and
\begin{equation}\label{cond_pen}
\sqrt{\sharp S_k(N)} (\widehat{\psi}(y,S_k(N))-\psi(y))=O_{\p}(1),\;\;N\to \infty.
\end{equation}
Let also \eqref{eta} hold.
Then, as $N\to \infty$,
\begin{equation*}
\sqrt{m_N}\Big(\dfrac{1}{m_N}\sum_{j=1}^{m_N}\dfrac{1}{\sqrt{K}}
\sum_{k=1}^K\sum_{i=0}^{2m-1}\sum_{i-m<|y|\le m}\Big[\widehat{\psi}(y,S_k(N))\ind\{A_N(i,j,k,y)\}\Big]
 \!-\! Err(f)\Big)
\stackrel{law}\longrightarrow Z_{0,\sigma^2}.
\end{equation*}
Here $A_N(i,j,k,y)=\{Y^{j+(k-1)q}=y, |f_{PA}(X^{j+(k-1)q})-y|>i\}$,
$Z_{0,\sigma^2}\sim \mathcal{N}(0,\sigma^2)$ and $\sigma^2$
was introduced in \eqref{sigma^2}.
\end{Th}

{\it Proof.} One can show that
$$
\sqrt{m_N}\big(\widehat{L}_{q}-Err(f)\big)
$$
$$
-\,\dfrac{\sqrt{m_N}}{\sqrt{K}}\sum_{k=1}^K\sum_{i=0}^{2m-1}\sum_{i-m<|y|\le m} \Big[
(\widehat\psi(y,{S_k(N)})-\psi(y))\p(Y=y,|f(X)-y|>i)
$$
$$
+\,\psi(y)\Big(\dfrac{1}{\sharp S_k(N)}\sum_{j\in S_k(N)}
H_N(y,i,j)
\Big)\Big]
\xrightarrow{{\sf P}}0
$$
as $N\to\infty$. Here $H_N(y,i,j)=\ind\{Y^j=y,|f(X^j)-y|>i\}-\p(Y=y,|f(X)-y|>i)$. For any
$i\in\{0,\ldots,2m-1\}$ and $y\in\mathbb Y$, the CLT for arrays of i.i.d. random variables with finite second moment implies that
\begin{equation*}
\dfrac{1}{\sqrt{\sharp S_k(N)}}\sum_{j\in S_k(N)}H_N(y,i,j)=O_{\p}(1),\;\;\;N\to\infty.
\end{equation*}
Since $m_N\slash \sharp S_k(N)\to0$ we get
$$
\dfrac{\sqrt{m_N}}{\sqrt{K}}\sum_{k=1}^K\sum_{i=0}^{2m-1}\sum_{i-m<|y|\le m} \psi(y)
\dfrac{1}{\sharp S_k(N)}\sum_{j\in S_k(N)}H_N(y,i,j)\xrightarrow{{\sf P}}0,\;\;\;N\to\infty.
$$
In a similar way in view of \eqref{cond_pen} one has
\begin{equation*}
\dfrac{\sqrt{m_N}}{\sqrt{K}}\sum_{k=1}^K\sum_{i=0}^{2m-1}\sum_{i-m<|y|\le m} \big
(\widehat\psi(y,{S_k(N)})-\psi(y)\big)\p(Y=y,|f(X)-y|>i)\xrightarrow{{\sf P}}0,\;\;\;N\to\infty.
\end{equation*}
Thus under conditions of Theorem \ref{Th2} the asymptotic behavior
of $\sqrt{m_N}\big(\widehat{L}_{m_N}-\widehat{L}_{q}\big)$ is the
same as for $\sqrt{m_N}\big(\widehat{L}_{m_N}-Err(f)\big)$.
$\square$

In Velez et al. (2007) the following choice of the penalty function $\psi$ was
proposed
\begin{equation*}
\psi(y)=\dfrac{c}{\p (Y=y)},\;\;y\in\mathbb Y,\;\;c=\text{const}>0.
\end{equation*}
This choice was justified in Bulinski (2012) for binary response $Y$. We will employ this penalty
function for nonbinary response as well, i.e. when $\mathbb
Y=\{-m,\ldots,0,\ldots,m\}$. Futher we assume that $\p (Y=y)>0$ for
all $y\in\mathbb Y$ and w.l.g. $c=1$.

Introduce $A_N(y,S_k(N))=\{Y^j\ne y,\; j\in S_k(N)\},\;N\in\mathbb
N,\;k\in\{1,\ldots,K\},\;y\in\mathbb Y$ and set (as usual $0/0:=0$)
\begin{equation}\label{pen_func}
\widehat{\p}_{S_k(N)}(Y=y):=\dfrac{\sum_{j\in S_k(N)}\ind\{Y^j=y\}}{\sharp
S_k(N)},\;\;\widehat\psi(y,S_k(N)):=\dfrac{\ind\{\Omega\setminus
A_N(y,S_k(N))\}}{\widehat{\p}_{S_k(N)}(Y=y)}.
\end{equation}

\begin{Cor}
The estimate $\widehat{\psi}$ defined by
\eqref{pen_func} satisfies
conditions of Theorem \ref{Th2}.
\end{Cor}

{\it Proof.} Fix arbitrary $y\in\mathbb{Y}$ and $k=1,\ldots,K$. One
can easily check that $\widehat{\psi}(y,S_k(N))$ is a strongly
consistent estimate of $\psi(y)$. Moreover, by CLT for arrays of
i.i.d. random variables we have
$$
\sqrt{\sharp\,S_k(N)}\big(\widehat{\p}_{S_k(N)}(Y=y)-\p(Y=y)\big)\xrightarrow{law}
Z_{0,\sigma_1^2(y)}\sim
\mathcal{N}\big(0,\sigma_1^2(y)\big),\;\;\;N\to\infty,
$$
where $\sigma_1^2(y)= \p(Y=y)(1-\p(Y=y))$. Taking into account that
$\widehat{\p}_{S_k(N)}(Y\!=\!y)\!\to\!\p(Y\!=\!y)$ a.s. and
$\sqrt{S_k(N)}\ind\{A_N(y, S_k(N))\}\xrightarrow{\sf P}0$ as $N\to
\infty$, one can write by Slutsky's lemma that
\begin{equation*}
\sqrt{\sharp S_k(N)}\big(\widehat{\psi}_{S_k(N)}(y)-\psi(y)\big)\xrightarrow{law}
Z_{0,\sigma_2^2(y)}\sim\mathcal{N}\big(0,\sigma_2^2(y)\big),\;\;\;N\to\infty,
\end{equation*}
where $\sigma_2^2(y)=(1-\p(Y=y))\p(Y=y)^{-3}$.
Thus \eqref{cond_pen} holds. Now we verify \eqref{eta}. Clearly,
$\widehat{\psi}(y,S_k(N))\le \sharp S_K(N)$ for any $N\in\mathbb N$.
Put $\varepsilon:=\min_{y\in\mathbb Y}\p(Y=y)$. Then
by the Hoeffding inequality
$$
\E |\widehat{\psi}(y,S_k(N))|^4=\E \Big[|\widehat{\psi}_{N,k}(y)|^4\ind\Big\{\big|\widehat{\p}_{S_k(N)}(Y=y)-\p(Y=y)\big|>\varepsilon\slash2\Big\}\Big]
$$
$$
+\E \Big[(\widehat{\psi}(y,S_k(N)))^4
\ind\Big\{\big|\widehat{\p}_{S_k(N)}(Y=y)-\p(Y=y)\big|\le\varepsilon\slash2\Big\}\Big]
$$
$$
\le 2(\sharp S_K(N))^4\exp\{-\sharp S_1(N)\varepsilon^2\slash2\}+2^4\slash\varepsilon^4.
$$
Thus we come to \eqref{eta}. $\square$

To simplify notation we will write in the following theorem $\widehat{Err}_K(f_{PA},\xi_N)$
for random variable introduced in \eqref{Errestmult}
replacing $\widehat{\psi}(y,S_k(N))$ by $\widehat{\psi}(y,\overline{S_k(N)})$, $y\in\mathbb{Y}$, $k=1,\ldots,K$. After such replacement in \eqref{xi_N} -- \eqref{sum_main} we
obtain the new row-wise exchangeable array $\{X_{N,j}, 1\leq j\leq q, N\in\mathbb{N}\}$
and therefore all established results hold true in this case.

\begin{Th}\label{Th3} Let $\varepsilon_N\to 0$ and $N^{1/2}\varepsilon_N \to \infty$ as $N\to \infty$.
Then, for each $K\in \mathbb{N}$, any vector $\beta=(k_1,\ldots k_r)$ with
$1\leq k_1<\ldots < k_r\leq n$, the corresponding
function $f=f^\beta$ and prediction algorithm defined by
$f_{PA}=\widehat{f}_{PA,\varepsilon}^{\beta}$, the following
relation holds:
\begin{equation}\label{CLT}
\sqrt{N}(\widehat{Err}_K(f_{PA},\xi_N) - Err(f))
\stackrel{law}\longrightarrow Z_{0,\sigma^2}\sim \mathcal{N}(0,\sigma^2),\;\;N\to \infty.
\end{equation}
Here $\sigma^2$ is  variance of the random variable
\begin{equation}\label{abc}
V= \sum_{i=0}^{2m-1}\sum_{i-m< |y|\le m}
\frac{\mathbb{I}\{Y=y\}}{{\sf P}(Y=y)} \left(\mathbb{I}\{|f(X)-
y|>i\} - {\sf P}(|f(X)- y|>i\big|Y=y)\right).
\end{equation}
\end{Th}

{\it Proof.} Set, for $f:\mathbb X\to\mathbb Y$ and $N\in \mathbb N$,
\begin{equation*}
T_N(f):=\dfrac{1}{K}\sum_{k=1}^K\dfrac{1}{\sharp S_k(N)}\sum_{i=1}^{2m-1}\sum_{i-m<|y|\le m}\psi(y)\sum_{j\in S_k(N)}\ind\{Y^j=y,|f(X^j)-y|>i\}.
\end{equation*}
The Slutsky lemma shows
that the limit behavior of the random
variables introduced in \eqref{eq_6} will be the same as for
random variables
$$
\rho_N:=\sqrt{N}(T_N(f)-Err(f))
$$
$$-\, \frac{\sqrt{N}}{K}\sum_{k=1}^K
\sum_{i=0}^{2m-1}\sum_{i-m< |y|\le m} \frac{(\widehat{\sf
P}_{\overline{S_k(N)}}(Y=y)- {\sf P}(Y=y)){\sf P}(Y=y,|f(X)- y|>i)}{{\sf
P}(Y=y)^2}
$$
$$
=\frac{\sqrt{N}}{K}\sum_{k=1}^K
\sum_{i=0}^{2m-1}\sum_{i-m< |y|\le m} \frac{1}{\sharp\,S_k(N)}
\sum_{j\in S_k(N)}\frac{\mathbb{I}\{Y^j=y,|f(X^j)-y|>i)}{{\sf P}(Y=y)}
$$
$$
- \, \frac{\sqrt{N}}{K}\sum_{k=1}^K
\sum_{i=0}^{2m-1}\sum_{i-m< |y|\le m} \frac{1}{\sharp\,\overline{S_k(N)}}
\sum_{j\in \overline{S_k(N)}}\mathbb{I}\{Y^j=y\}\frac{{\sf P}(Y=y, |f(X)-y|>i\}}{({\sf P}(Y=y))^2}.
$$
Let $a_k, b_k$, $k=1,\ldots,K$, be any real numbers. We use the following simple observation
$$
\frac{1}{\sharp\,S_k(N)}\sum_{k=1}^K a_k + \frac{1}{\sharp\,\overline{S_k(N)}}\sum_{l=1,\ldots,K;l\neq k}b_l
=\sum_{k=1}^K\left(\frac{a_k}{\sharp\,S_k(N)} + b_k\sum_{l=1,\ldots,K;l\neq k}
\frac{1}{\sharp\,\overline{S_l(N)}}\right).
$$
Combining the latter formulas we can write
$$
\rho_N=\frac{\sqrt{N}}{K}\sum_{k=1}^K \sum_{j\in
S_k(N)}\left(\frac{V^j_1}{\sharp\,S_k(N)}+ V^j_2\sum_{l=1,\ldots,K;l\neq k}
\frac{1}{\sharp\,\overline{S_l(N)}}\right)
$$
where
$$
V^j_1=\sum_{i=0}^{2m-1}\sum_{i-m< |y|\le m}\frac{\mathbb{I}\{Y^j=y, |f(X^j)-y|>i\}}{{\sf P}(Y=y)},
$$
$$
V^j_2=-\sum_{i=0}^{2m-1}\sum_{i-m< |y|\le m}
\frac{\mathbb{I}\{Y^j=y\}{\sf P}(Y=y,|f(X)-y|>i)}{({\sf P}(Y=y))^2}.
$$
Take any $k\in\{1,\ldots,K\}$ and employ CLT for an array of bounded centered i.i.d. random variables
$\{V^j_1 + V^j_2, j\in S_k(N), N\in\mathbb{N}\}$. Then
$$
\frac{1}{\sqrt{\sharp\,S_k(N)}}\sum_{j\in S_k(N)}(V^j_1 + V^j_2) \stackrel{law}\longrightarrow
Z_{0,\sigma^2}\sim \mathcal{N}(0,\sigma^2),\;\;N\to \infty,
$$
where $\sigma^2= {\sf var} (V_1^j+V_2^j)$. Note now that, for each $k=1,\ldots,K$,
$$
\frac{N}{\sharp\,S_k(N)}\to \frac{1}{K},\;\;\; \sum_{l=1,\ldots,K;l\neq k}
\frac{N}{\sharp\,\overline{S_l(N)}}\to \frac{1}{K},\;\;\;N\to \infty.
$$
For each $N\in\mathbb{N}$, the
families of random variables $\{V^j_1 + V^j_2, j\in S_k(N)\}$, $k=1,\ldots,K$, are
independent. Thus we come to the following relation
$$
\rho_N \stackrel{law}\longrightarrow
Z_{0,\sigma^2}\sim \mathcal{N}(0,\sigma^2),\;\;N\to \infty.
$$
Obviously we can write $\sigma^2= {\sf var} \,V$ where
$V$ is introduced in \eqref{abc}.
The proof is complete. $\square$

\begin{Rem}\label{Rem4}
It is not difficult  to construct the consistent estimates
$\widehat{\sigma}_N$ of unknown $\sigma$ appearing in \eqref{CLT}.
Therefore (if $\sigma^2\neq 0$) we can claim that under conditions
of Theorem~\ref{Th3}
$$
\frac{\sqrt{N}}{\widehat{\sigma}_N}(\widehat{E}rr_K(f_{PA},\xi_N) -
Err(f)) \stackrel{law}\longrightarrow \frac{Z}{\sigma}\sim
\mathcal{N}(0,1),\;\;N\to \infty.
$$
\end{Rem}

\begin{Rem}\label{Rem5}
It is interesting to compare Theorem \ref{Th3}
with simulations corresponding to Example~3 when $N$ is rather small (e.g., $N=200$).
In this case the choice of $\widehat{\psi}(y,\overline{S_k(N)})$
can lead to better identification of significant factors.
\end{Rem}

\vskip1cm

\noindent BIBLIOGRAPHY
\vskip 3mm

\noindent Arlot, S.,  Celisse,  A. (2010).  A survey of cross-validation
procedures for model selection. {\it Statist. Surv.}. \textbf{4},
40--79.
\vskip 3mm

\noindent Berti, Patrizia; Pratelli, Luca; Rigo, Pietro. (2004).
Limit theorems for a class of identically distributed random variables.
{\it Ann. Probab.} 32, no. 3, 2029--2052.
\vskip 3mm

\noindent Billingsley P. (1968). {\it Convergence of Probability Measures.} John Wiley
$\&$ Sons Inc., New York.
\vskip 3mm

\noindent Bulinski, A.V. (2014). On foundation of the dimensionality reduction method
for explanatory variables. {\it J. Math. Sci.}, DOI 10.1007/s10958-014-1838-7.
\vskip 3mm

\noindent Bulinski, A.V. (to appear, 2014). Central limit theorem related to MDR method.:
Proceedings of the Fields Institute International Symposium on Asymptotic Methods in Stochastics,
in Honour of Mikl\'{o}s Cs\"{o}rg\H{o}'s Work
on the occasion of his anniversary.
arXiv:1301.6609 [math.PR].
\vskip 3mm

\noindent Bulinski, A., Butkovsky, O., Sadovnichy, V.,
Shashkin, A., Yaskov, P., Balatskiy, A., Samokhodskaya, L., Tkachuk,
V. (2012). Statistical methods of SNP data analysis and applications. {\it Open
Journal of Statistics.} \textbf{2}(1), 73--87.
\vskip 3mm

\noindent Bulinski A.V., Rakitko A.S. (2014). Estimation of nonbinary random  response.
{\it Dokl. Math.}, \textbf{455}(6), 1--5.
\vskip 3mm

\noindent Chernoff, H., Teicher, H. (1958).
A Central Limit Theorem for Sums of Interchangeable Random Variables.
{\it Ann. Math. Statist.} \textbf{29}(1),118--130.
\vskip 3mm

\noindent Golland, P., Liang, F., Mukherjee, S., Panchenko, D. (2005). {\it Permutation Tests for Classification.} LNCS, \textbf{3559}, 501--515.
\vskip 3mm

\noindent Hastie T., Tibshirani R. and Friedman J. (2008). {\it The Elements of Statistical Learning; Data Mining,
Inference and Prediction.} Springer, New York. Second edition.
\vskip 3mm

\noindent Lee S., Epstein M.P., Duncan R. and Lin X. (2012).
Sparse Principal Component Analysis for Identifying Ancestry-Informative
Markers in Genome Wide Association Studies. {\it Genet. Epidemiol.}
\textbf{36}, 293--302.
\vskip 3mm

\noindent Moore J.B., , Asselbergs F.W.
and Williams S.M. (2010). Bioinformatics challenges for genome-wide association studies.
{\it Bioinformatics.} \textbf{26}, 445--455.
\vskip 3mm

\noindent Ritchie, M.D., Hahn, L.W., Roodi, N., Bailey, R.L.,
Dupont, W.D., Parl, F.F., Moore, J.H. (2001). Multifactor-dimensionality
reduction reveals high-order interactions among estrogen-metabolism
genes in sporadic breast cancer.  {\it Am. J. Hum Genet.} \textbf{69}(1),
138--147.
\vskip 3mm

\noindent R{\"o}llin, A. (2013). Stein's method in high dimensions with applications.
{\it Ann. Inst. Henri Poincar\'{e} Probab. Stat.}
\textbf{49}(2), 529--549.
\vskip 3mm

\noindent Schwender, H., Ruczinski, I. (2010). Logic regression and its
extensions. {\it Adv. Genet.}. \textbf{72}, 25-45.
\vskip 3mm

\noindent Sikorska K., Lesaffre E., Groenen P.F.G.,
and Eilers P.H.C. (2013).
GWAS on your notebook: fast semi-parallel
linear and logistic regression for genome-wide
association studies. {\it BMC Bioinformatics}, \textbf{14}:166.
\vskip 3mm

\noindent Tibshirani R.J. and Taylor J. (2012).
Degrees of freedom in lasso problems. {\it Ann. Statist.}
\textbf{40}, 1198--1232.
\vskip 3mm

\noindent Velez, D.R., White, B.C., Motsinger, A.A., Bush, W.S.,
Ritchie, M.D., Williams, S.M., Moore, J.H. (2007). A balanced accuracy
function for epistasis modeling in imbalanced datasets using
multifactor dimensionality reduction. {\it Genet. Epidemiol.}.
\textbf{31}(4), 306--315.
\vskip 3mm

\noindent Visscher, P.M., Brown, M.A., McCarthy, M.I., Yang, J. (2012).
Five Years of GWAS Discovery. {\it  Am. J. Hum. Genet.}. \textbf{90}, 7--24.
\vskip 3mm

\noindent Weber N.C. (1980). A martingale approach to central limit theorems for exchangeable random variables. {\it J. Appl. Probab.}, \textbf{17}(3), 662-673.
\vskip 3mm

\end{document}